\documentclass[12pt]{article}
\usepackage{latexsym,amssymb,amsmath,enumerate,float,geometry,cite}
\newtheorem{theorem}{Theorem}

\newtheorem{lemma}[theorem]{Lemma}
\newtheorem{conjecture}[theorem]{Conjecture}

\newtheorem{claim}{Claim}

\begin{document}

\title{\Large Exponential Domination in Subcubic Graphs}
\author{St\'{e}phane Bessy$^1$, Pascal Ochem$^1$, and Dieter Rautenbach$^2$}
\date{}
\maketitle
\vspace{-10mm}
\begin{center}
{\small
$^1$ 
Laboratoire d'Informatique, de Robotique et de Micro\'{e}lectronique de Montpellier (LIRMM),\\
Montpellier, France, \texttt{stephane.bessy@lirmm.fr,pascal.ochem@lirmm.fr}\\[3mm]
$^2$ Institute of Optimization and Operations Research, Ulm University,\\
Ulm, Germany, \texttt{dieter.rautenbach@uni-ulm.de}}
\end{center}

\begin{abstract}
As a natural variant of domination in graphs, 
Dankelmann et al.~[Domination with exponential decay, Discrete Math. 309 (2009) 5877-5883] introduce exponential domination,
where vertices are considered to have some dominating power that decreases exponentially with the distance,
and the dominated vertices have to accumulate a sufficient amount of this power emanating from the dominating vertices.
More precisely, if $S$ is a set of vertices of a graph $G$, 
then $S$ is an exponential dominating set of $G$ if 
$\sum\limits_{v\in S}\left(\frac{1}{2}\right)^{{\rm dist}_{(G,S)}(u,v)-1}\geq 1$
for every vertex $u$ in $V(G)\setminus S$,
where ${\rm dist}_{(G,S)}(u,v)$ is the distance between $u\in V(G)\setminus S$ and $v\in S$ in the graph $G-(S\setminus \{ v\})$.
The exponential domination number $\gamma_e(G)$ of $G$ is the minimum order of an exponential dominating set of $G$.

In the present paper we study exponential domination in subcubic graphs.
Our results are as follows:
If $G$ is a connected subcubic graph of order $n(G)$, then
$$\frac{n(G)}{6\log_2(n(G)+2)+4}\leq \gamma_e(G)\leq \frac{1}{3}(n(G)+2).$$
For every $\epsilon>0$, there is some $g$ such that $\gamma_e(G)\leq \epsilon n(G)$ for every cubic graph $G$ of girth at least $g$.
For every $0<\alpha<\frac{2}{3\ln(2)}$, there are infinitely many cubic graphs $G$ with 
$\gamma_e(G)\leq \frac{3n(G)}{\ln(n(G))^{\alpha}}$.
If $T$ is a subcubic tree, then $\gamma_e(T)\geq \frac{1}{6}(n(T)+2).$
For a given subcubic tree, $\gamma_e(T)$ can be determined in polynomial time.
The minimum exponential dominating set problem is APX-hard for subcubic graphs.
\end{abstract}

{\small
\begin{tabular}{lp{12.5cm}}
\textbf{Keywords:} & domination, exponential domination, subcubic graph, cubic graph, girth, cage\\
\textbf{MSC2010:} & 05C69
\end{tabular}
}

\pagebreak
 
\section{Introduction}

We consider finite, simple, and undirected graphs, and use standard notation and terminology.

A set $D$ of vertices of a graph $G$ is {\it dominating} if every vertex not in $D$ has a neighbor in $D$.
The {\it domination number} $\gamma(G)$ of $G$, defined as the minimum cardinality of a dominating set, 
is one of the most well studied quantities in graph theory \cite{hhs}.
As a natural variant of this classical notion, Dankelmann et al.~\cite{ddems} introduce exponential domination,
where vertices are considered to have some dominating power that decreases exponentially by the factor $\frac{1}{2}$ with the distance,
and the dominated vertices have to accumulate a sufficient amount of this power emanating from the dominating vertices.
As a motivation of their model they mention information dissemination within social networks, where
the impact of information decreases every time it is passed on.

Before giving the precise definitions for exponential domination, 
we point out that it shares features with several other well studied domination notions,
such as, for example,
{\it $k$-domination} \cite{cr,cgs,dghpv,fhv,fj,ha,hv,rv}, where several vertices contribute to the domination of an individual vertex,
{\it distance-$k$-domination} \cite{adr,bz,chmm,hmv,he,tx}, where vertices dominate others over some distance, and 
{\it broadcast domination} \cite{chm,dehhh,e,hl}, 
where different dominating vertices contribute differently to the domination of an individual vertex depending on the relevant distances.

Let $G$ be a graph.
The vertex set and the edge set of $G$ are denoted by $V(G)$ and $E(G)$, respectively.
The order $n(G)$ of $G$ is the number of vertices of $G$,
and the size $m(G)$ of $G$ is the number of edges of $G$.
For two vertices $u$ and $v$ of $G$, let ${\rm dist}_G(u,v)$ be the distance in $G$ between $u$ and $v$,
which is the minimum number of edges of a path in $G$ between $u$ and $v$.
If no such path exists, then let ${\rm dist}_G(u,v)=\infty$.
An endvertex is a vertex of degree at most $1$.
For a rooted tree $T$, and a vertex $u$ of $T$, 
let $T_u$ denote the subtree of $T$ rooted in $u$ that contains $u$ as well as all descendants of $u$.
A leaf of a rooted tree is a vertex with no children.
For non-negative integers $d_0,d_1,\ldots,d_k$, 
let $T(d_0,d_1,\ldots,d_k)$ be the rooted tree of depth $k+1$
in which all vertices at distance $i$ from the root have exactly $d_i$ children for every $i$ with $0\leq i\leq k$.
A rooted tree is binary if every vertex has at most two children, 
and a binary tree is full if every vertex other than the leaves has exactly two children.
For a positive integer $k$, let $[k]$ be the set of positive integers at most $k$.

Let $S$ be a set of vertices of $G$.
For two vertices $u$ and $v$ of $G$ with $u\in S$ or $v\in S$,
let ${\rm dist}_{(G,S)}(u,v)$ be the minimum number of edges of a path $P$ in $G$ between $u$ and $v$ 
such that $S$ contains exactly one endvertex of $P$ and no internal vertex of $P$.
If no such path exists, then let ${\rm dist}_{(G,S)}(u,v)=\infty$.
Note that, if $u$ and $v$ are distinct vertices in $S$, then ${\rm dist}_{(G,S)}(u,u)=0$ and ${\rm dist}_{(G,S)}(u,v)=\infty$.

For a vertex $u$ of $G$, let 
$$w_{(G,S)}(u)=\sum\limits_{v\in S}\left(\frac{1}{2}\right)^{{\rm dist}_{(G,S)}(u,v)-1},$$
where $\left(\frac{1}{2}\right)^{\infty}=0$.
Note that $w_{(G,S)}(u)=2$ for $u\in S$.

If $w_{(G,S)}(u)\geq 1$ for every vertex $u$ of $G$, then $S$ is an {\it exponential dominating set} of $G$.
The {\it exponential domination number $\gamma_e(G)$} is the minimum order of an exponential dominating set of $G$,
and an exponential dominating set of $G$ of order $\gamma_e(G)$ is {\it minimum}.
By definition, every dominating set is also an exponential dominating set, 
which implies $\gamma_e(G)\leq \gamma(G)$ for every graph $G$.
Dankelmann et al.~\cite{ddems} also consider a {\it porous} version, 
where the term ``$\left(\frac{1}{2}\right)^{{\rm dist}_{(G,S)}(u,v)-1}$'' in the definition of $w_{(G,s)}(u)$
is replaced by ``$\left(\frac{1}{2}\right)^{{\rm dist}_{G}(u,v)-1}$''.
Note that in the porous version, the different vertices in $S$ do not block each others influence.

\medskip

\noindent In the present paper we focus on exponential domination in subcubic graphs,
which is a special case that displays several interesting features.
The intuitive reason for this is that, by definition, the dominating power halves with every additional unit of distance, 
while, at least in subcubic graphs, 
the number of vertices at a certain distance from a given vertex at most doubles with every additional unit of distance,
that is, the product of these two quantities is bounded.

The following lemma makes this vague observation more precise.

\begin{lemma}\label{lemmad31}
Let $G$ be a graph of maximum degree at most $3$,
and let $S$ be a set of vertices of $G$.

If $u$ is a vertex of degree at most $2$ in $G$, then $w_{(G,S)}(u)\leq 2$ with equality if and only if 
$u$ is contained in a subgraph $T$ of $G$ that is a tree, 
such that rooting $T$ in $u$ yields a full binary tree and $S\cap V(T)$ is exactly the set of leaves of $T$.
\end{lemma}
For a general graph $G$, the value of $w_{(G,S)}(u)$ is not bounded from above, 
and may, in particular, be arbitrarily large because of vertices in $S$ at any distance from $u$.
This turns exponential domination into a non-local problem.
In contrast to that, 
Lemma \ref{lemmad31} implies $w_{(G,S)}(u)\leq 3$ for every subcubic graph $G$, every set $S$ of its vertices, and every vertex $u$ of $G$,
that is, the accumulated effect of arbitrarily many vertices at any distance from $u$ 
can not be substantially larger than the effect of a single neighbor of $u$ in $S$,
which already implies $w_{(G,S)}(u)\geq 1$.
This somewhat localizes exponential domination in subcubic graphs, 
which makes it more tractable in many aspects.

Let $S$ be a set of vertices of a graph $G$,
and let $u$ and $v$ be distinct vertices such that $u\not\in S$ and $w_{(G,S)}(v)\geq 1$.
If $G$ is subcubic, then Lemma \ref{lemmad31} implies $w_{(G,S\cup \{ u\})}(v)\geq 1$,
which does not hold in general.
In particular, a superset of an exponential dominating set of a not necessarily subcubic graph may not be an exponential dominating set.
If $w_{(G,S)}(u)\leq 2$ though, then $w_{(G,S\cup \{ u\})}(v)\geq 1$ follows also for not necessarily subcubic graphs. 

The main results of Dankelmann et al.~\cite{ddems} are as follows.

\begin{theorem}[Dankelmann et al.~\cite{ddems}]\label{theoremd}
If $G$ is a connected graph of diameter ${\rm diam}(G)$, then 
$$\frac{1}{4}({\rm diam}(G)+2)\leq \gamma_e(G)\leq \frac{2}{5}(n(G)+2).$$
\end{theorem}
While the lower bound is tight for paths of order $2$ mod $4$, 
a tight upper bound for connected graphs in terms of their order is still unknown.
In \cite{bor} we provide further bounds, and, in particular, strengthen the upper bound to 
$\gamma_e(G)\leq \frac{43}{108}(n(G)+2)$, which is still surely not best possible.
For subcubic graphs though, we obtain a tight upper bound.

\begin{theorem}\label{theoremub3}
If $G$ is a connected graph of maximum degree at most $3$, then $\gamma_e(G)\leq \frac{1}{3}(n(G)+2)$.
\end{theorem}
It seems even possible to characterize the extremal graphs for Theorem \ref{theoremub3},
and we formulate an explicit conjecture at least for the extremal trees.

Since the diameter of a subcubic graph $G$ is at least linear in $\log (n(G))$,
Theorem \ref{theoremd} implies $\gamma_e(G)\geq \Omega(\log (n(G)))$.
We improve this lower bound as follows.

\begin{theorem}\label{theoremlb3}
If $T$ is a tree of maximum degree at most $3$, then $\gamma_e(T)\geq \frac{1}{6}(n(T)+2)$.
\end{theorem}

\begin{theorem}\label{theoremlb3gen}
If $G$ is a graph of maximum degree at most $3$, then $\gamma_e(G)\geq \frac{n(G)}{6\log_2(n(G)+2)+4}$.
\end{theorem}
Our next result implies that Theorem \ref{theoremlb3gen} is not far from being best possible.

\begin{theorem}\label{theoremgirth}
For every $\epsilon>0$, there is some $g$ such that $\gamma_e(G)\leq \epsilon n(G)$ for every cubic graph $G$ of girth at least $g$.
Furthermore, for every $\alpha$ with $0<\alpha<\frac{2}{3\ln(2)}$, 
there are infinitely many cubic graphs $G$ with 
$\gamma_e(G)\leq \frac{3n(G)}{\ln(n(G))^{\alpha}}$.
\end{theorem}
Imposing a stronger condition actually allows to derive a linear lower bound.

\begin{theorem}\label{theoremlb3sp}
Let $G$ be a graph of order at least $3$ and maximum degree at most $3$.
If $S$ is a set of vertices of $G$ such that $w_{(G,S)}(u)\geq 3$ for every $u\in V(G)\setminus S$,
then $|S|\geq \frac{1}{4}(n(G)+6)$.
\end{theorem}
While Dankelmann et al.~\cite{ddems} do not comment on the complexity of the exponential domination number in general,
they explicitly ask whether there is a polynomial time algorithm that computes the exponential domination number of a given tree.
Relying on Lemma \ref{lemmad31}, we obtain such an algorithm for subcubic trees.

\begin{theorem}\label{theoremalg3}
For a given tree $T$ of maximum degree at most $3$, 
$\gamma_e(T)$ can be determined in polynomial time.
\end{theorem}
Finally, we establish a hardness result.

\begin{theorem}\label{theoremnpc}
The problem to determine a minimum exponential dominating set of a given subcubic graph is APX-hard.
\end{theorem}
All proofs and further discussion are in the next section.
A third section summarizes several open problems related to our results.

\section{Proofs of the results}

{\it Proof of Lemma \ref{lemmad31}:} For a non-negative integer $i$, let 
$$V_i=\{ v\in V(G):{\rm dist}_{(G,S\cup \{ v\})}(u,v)=i\},$$ 
and let $n_i=|V_i|$ and $s_i=|S\cap V_i|$.
Note that for every vertex $v$ in $V_i$, there is a path $P$ in $G$ between $u$ and $v$
such that $V(P)\setminus \{ v\}$ does not intersect $S$,
and the minimum length of such a path is exactly $i$.
Trivially, $V_0=\{ u\}$ and thus $n_0=1$.
Since every vertex in $V_{i+1}$ has a neighbor in $V_i\setminus S$, 
the degree conditions imply $n_{i+1}\leq 2(n_i-s_i)$.

\begin{claim}\label{claim1}
For every non-negative integer $k$, 
\begin{eqnarray}\label{e1}
n_k & \leq & 2^k\left(1-\sum\limits_{i=0}^{k-1}\frac{s_i}{2^i}\right)
\end{eqnarray}
with equality if and only if $n_i=2(n_{i-1}-s_{i-1})$ for every $i\in [k]$.
\end{claim}
{\it Proof of Claim \ref{claim1}:} We prove the claim by induction on $k$.
Since $n_0=1$, the claim holds for $k=0$. 
For $k>0$, we obtain, by induction,
\begin{eqnarray*}
n_k & \leq & 2(n_{k-1}-s_{k-1})\\
& \leq & 2\left(2^{k-1}\left(1-\sum\limits_{i=0}^{k-2}\frac{s_i}{2^i}\right)-s_{k-1}\right)\\
& = & 2^k\left(1-\sum\limits_{i=0}^{k-1}\frac{s_i}{2^i}\right),
\end{eqnarray*}
and (\ref{e1}) follows.
Furthermore, we have equality in (\ref{e1}) if and only if 
$n_k=2(n_{k-1}-s_{k-1})$
and
$n_{k-1}=2^{k-1}\left(1-\sum\limits_{i=0}^{k-2}\frac{s_i}{2^i}\right)$.
By induction, this is equivalent to $n_i=2(n_{i-1}-s_{i-1})$ for every $i\in [k]$.
$\Box$

\medskip

\noindent Since $s_k\leq n_k$ for every non-negative integer $k$, Claim \ref{claim1} implies $\sum\limits_{i=0}^{k}\frac{s_i}{2^{i-1}}\leq 2$ and then
$$w_{(G,S)}(u)=\sum\limits_{i=0}^{\infty}\frac{s_i}{2^{i-1}}\leq 2.$$
Furthermore, if $w_{(G,S)}(u)=2$, then, since $G$ is finite, there is some non-negative integer $k$ such that  $\sum\limits_{i=0}^{k}\frac{s_i}{2^{i-1}}=2$, that is
$$s_k=2^k\left(1-\sum\limits_{i=0}^{k-1}\frac{s_i}{2^i}\right)\le n_k$$
which, by Claim \ref{claim1}, implies
$s_k=n_k$ and $n_i=2(n_{i-1}-s_{i-1})$ for every $i\in [k]$.
Since $u$ has degree at most $2$ and $G$ has maximum degree at most $3$, this implies the existence of the desired subtree $T$.
Conversely, if a tree $T$ as in the statement exists, then, clearly, $w_{(G,S)}(u)=2$. $\Box$

\medskip

\noindent Before we proceed to the proof of Theorem \ref{theoremub3}, 
we establish a lemma concerning specific reductions.

\begin{lemma}\label{lemmad32}
Let $G$ be a graph of maximum degree at most $3$.
Let $u$ be a vertex of $G$, and for a positive integer $k\in [2]$, let $v_1,\ldots,v_k$ be neighbors of $u$ in $G$
such that, for every $i\in [k]$, the component $T_i$ of $G-u$ that contains $v_i$ is a tree.


\begin{enumerate}[(i)]
\item If $k=2$, and $v_1$ and $v_2$ have degree $1$ in $G$, then $\gamma_e(G)=\gamma_e(G-v_2)$.
\item If $k=2$, $v_1$ has degree $1$ in $G$, and $T_2$ has depth $1$, then $\gamma_e(G)=\gamma_e(G-(\{ v_1\}\cup V(T_2)))+1$.
\item If $k=1$, $v_1$ has degree $3$ in $G$, $T_1$ has depth 2, and the two children $w_1$ and $w_2$ of $v_1$ in $T_1$ are both no leaves of $T_1$, 
then $\gamma_e(G)=\gamma_e(G-(V(T_1)\setminus \{ v_1,w_1,w_2\}))+1=\gamma_e(G-(V(T_1)\setminus \{ v_1,w_1\}))+1$.
\item If $k=1$ and $T_1\cong T(1,1)$, 
then $\gamma_e(G)\leq \gamma_e(G-V(T_1))+1$.
\end{enumerate}
\end{lemma}
{\it Proof:} (i) Note that an exponential dominating set of $G$ that contains no endvertices
is also an exponential dominating set of $G-\{ v_2\}$, and vice versa.
Since both $G$ as well as $G-\{ v_2\}$ have minimum exponential dominating sets that contain no endvertices, (i) follows.

\medskip

\noindent (ii) Let $G'=G-(\{ v_1\}\cup V(T_2))$.
Lemma \ref{lemmad31} implies that $G$ has a minimum exponential dominating set $S$ that contains $v_2$ 
but no other vertex from $\{ v_1\}\cup V(T_2)$. 
Let $S'=S\setminus \{ v_2\}$.
If $u\in S'$, then $S'$ is an exponential dominating set of $G'$.
Now, let $u\not\in S'$.
Since $1\leq w_{(G,S)}(v_1)=\frac{1}{2}w_{(G',S')}(u)+\frac{1}{2}$, we obtain $w_{(G',S')}(u)\geq 1$.
This implies that $u$ has a neighbor $x$ distinct from $v_1$ and $v_2$, and that $w_{(G',S')}(x)\geq 2$.
By Lemma \ref{lemmad31}, $G-u$ contains a subgraph $T$ that is a tree,
such that rooting $T$ in $x$ yields a full binary tree and $S'\cap V(T)$ is exactly the set of leaves of $T$.
Since $G$ is subcubic, this implies that, also in this case, $S'$ is an exponential dominating set of $G'$, and hence, 
$\gamma_e(G)\geq\gamma_e(G')+1$.
Conversely, if $S'$ is a minimum exponential dominating set of $G'$, 
then $S'\cup \{ v_2\}$ is an exponential dominating set of $G$, which implies $\gamma_e(G)\leq\gamma_e(G')+1$.

\medskip

\noindent (iii) Let $G'=G-(V(T_1)\setminus \{ v_1,w_1,w_2\})$.
Lemma \ref{lemmad31} implies that $G$ has a minimum exponential dominating set $S$ that contains $w_1$ and $w_2$
but no other vertex from $V(T_1)$.
Since $(S\setminus \{ w_1,w_2\})\cup \{ v_1\}$ is an exponential dominating set of $G'$, 
we have $\gamma_e(G)\geq\gamma_e(G')+1$.
Conversely, 
Lemma \ref{lemmad31} implies that $G'$ has a minimum exponential dominating set $S'$ that contains $v_1$
but neither $w_1$ nor $w_2$.
Since $(S'\setminus \{ v_1\})\cup \{ w_1,w_2\}$ is an exponential dominating set of $G$, 
we have $\gamma_e(G)\leq\gamma_e(G')+1$.
The equality $\gamma_e(G-(V(T_1)\setminus \{ v_1,w_1,w_2\}))=\gamma_e(G-(V(T_1)\setminus \{ v_1,w_1\}))$
follows from (i).

\medskip

\noindent (iv) Since adding the child of $v_1$ in $T_1$ to an exponential dominating set of $G-V(T_1)$ 
yields an exponential dominating set of $T$, (iv) follows. $\Box$

\medskip

\noindent {\it Proof of Theorem \ref{theoremub3}:}
If $H$ is a spanning subgraph of $G$, then $\gamma_e(G)\leq \gamma_e(H)$.
Therefore, it suffices to prove Theorem \ref{theoremub3} if $G$ is a tree $T$.
The proof is by induction on the order of $T$.
If $T$ has diameter at most $2$, then $\gamma_e(T)=1$, and the bound follows.
Hence, we may assume that the diameter of $T$ is at least $3$, and hence, $n(T)\geq 4$.
Root $T$ in a vertex of maximum eccentricity.
Note that $T$ has depth at least $3$.

If some vertex of $T$ has two children that are leaves,
then Lemma \ref{lemmad32}(i) implies the existence of a tree $T'$ with $n(T')<n(T)$ and $\gamma_e(T)=\gamma_e(T')$.
By induction, $\gamma_e(T)\leq \frac{1}{3}(n(T')+2)<\frac{1}{3}(n(T)+2)$.
Hence, we may assume that no vertex of $T$ has two children that are leaves.

Let $u$ be a vertex of $T$ such that $T_u$ has depth $2$. 

If $u$ has a child that is a leaf, then Lemma \ref{lemmad32}(ii) implies the existence of a tree $T'$ with $n(T')=n(T)-3$ and $\gamma_e(T)=\gamma_e(T')+1$.
By induction, $\gamma_e(T)\leq \frac{1}{3}(n(T')+2)+1=\frac{1}{3}(n(T)+2)$.
Hence, we may assume that no child of $u$ is a leaf.
If $u$ has two children, then Lemma \ref{lemmad32}(iii) implies the existence of a tree $T'$ with $n(T')=n(T)-3$ and $\gamma_e(T)=\gamma_e(T')+1$.
By induction, $\gamma_e(T)\leq \frac{1}{3}(n(T')+2)+1=\frac{1}{3}(n(T)+2)$.
Since $T$ has maximum degree at most $3$ and depth at least $3$, we may assume that $u$ has exactly one child,
which implies $T_u\cong T(1,1)$.
Now, by induction, Lemma \ref{lemmad32}(iv) implies,
$\gamma_e(T)\leq \gamma_e(T-V(T_u))+1\leq \frac{1}{3}(n(T')+2)+1=\frac{1}{3}(n(T)+2)$,
which completes the proof. $\Box$

\medskip

\noindent As stated in the introduction, it seems possible to characterize the extremal graphs for Theorem \ref{theoremub3}.
Most of this hope is based on the equalities in Lemma \ref{lemmad32}(i), (ii), and (iii),
which should allow to relate extremal graphs of different orders in a controlled way.
In order to phrase a precise conjecture at least for the extremal trees,
we define three operations. 

Let $T$ and $T'$ be two trees.
\begin{itemize} 
\item $T$ arises by applying {\it Operation 1 to $T'$ at vertex $u$}
if $T$ contains a path $v_1uv_2w$ such that $v_1$ and $w$ are endvertices of $T$, $v_2$ has degree $2$ in $T$, and $T'=T-\{ v_1,v_2,w\}$.
\item $T$ arises by applying {\it Operation 2 to $T'$ at vertex $u$}
if $T$ contains an edge $xu$ such that rooting $T$ in $x$, we have $T_u\cong T(2,1)$,
and $T'=T-(V(T_u)\setminus \{ u,v\})$, where $v$ is a child of $u$ in $T_u$.
\item $T$ arises by applying {\it Operation 3 to $T'$ at vertex $u$}
if $T$ contains an edge $uv$ such that rooting $T$ in $u$, we have $T_v\cong T(1,1)$,
and $T'=T-V(T_v)$.
\end{itemize}

\begin{conjecture}\label{conjectured3}
A tree $T$ of maximum degree at most $3$ satisfies $\gamma_e(T)=\frac{1}{3}(n(T)+2)$
if and only if it arises from $K_1$ by iteratively applying
Operation 1 at vertices of degree at most $1$,
Operation 2 at vertices of degree $2$, and
Operation 3 at vertices of degree $2$.
\end{conjecture}
Using Lemma \ref{lemmad32} it is not difficult to show that every extremal tree can be constructed as stated in the conjecture;
the hard part is to show the converse.
It also follows easily that all trees of maximum degree at most $3$ 
in which every vertex of degree $2$ is adjacent to an endvertex
are extremal.

\medskip

\noindent {\it Proof of Theorem \ref{theoremlb3}:}
The proof is by induction on the order of $T$.
Since the bound is trivial for trees of diameter at most $2$, 
we may assume that the diameter of $T$ is at least $3$, and hence, $n(T)\geq 4$.
Root $T$ in a vertex of maximum eccentricity.

Let $v$ be a vertex of $T$ such that $T_v$ has depth $2$.
If one child of $v$ is a leaf, then removing the at most four descendants of $v$ yields a tree $T'$   
such that, by Lemma \ref{lemmad32}(ii) and induction,
$\gamma_e(T)=\gamma_e(T')+1\geq \frac{1}{6}(n(T')+2)+1\geq \frac{1}{6}(n(T)-4+2)+1>\frac{1}{6}(n(T)+2)$.
Hence, we may assume that no child of $v$ is a leaf.
If $v$ has two children, say $w_1$ and $w_2$, then removing the at most four descendants of $w_1$ and $w_2$ yields a tree $T'$
such that, by Lemma \ref{lemmad32}(iii) and induction,
$\gamma_e(T)=\gamma_e(T')+1\geq \frac{1}{6}(n(T')+2)+1\geq \frac{1}{6}(n(T)-4+2)+1>\frac{1}{6}(n(T)+2)$.
Hence, we may assume that $v$ has degree $2$ in $T$.

Let $u$ be a vertex of $T$ such that $T_u$ has depth $3$.
Let $v$ be a child of $u$ such that $T_v$ has depth $2$. 
By the previous observations, $v$ has degree $2$ in $T$.
Let $w$ be the child of $v$ in $T$.
If $u$ has a child that is a leaf, 
then, by Lemma \ref{lemmad31}, some minimum exponential dominating set $S$ of $T$ contains $u$ and $w$.
Since $S\setminus \{ w\}$ is an exponential dominating set of $T'=T-V(T_v)$, we obtain, by induction, 
$\gamma_e(T)\geq \gamma_e(T')+1\geq \frac{1}{6}(n(T')+2)+1\geq \frac{1}{6}(n(T)-4+2)+1>\frac{1}{6}(n(T)+2)$.
If $u$ has a child $v'$ such that $T_{v'}$ has depth $1$,
then, by Lemma \ref{lemmad31}, some minimum exponential dominating set $S$ of $T$ contains $v'$ and $w$.
If $T'$ arises from $T$ by removing the at most five descendants of $v$ and $v'$, 
then $(S\setminus \{ v',w\})\cup \{ u\}$ is an exponential dominating set of $T'$, and we obtain, by induction, 
$\gamma_e(T)\geq \gamma_e(T')+1\geq \frac{1}{6}(n(T')+2)+1\geq \frac{1}{6}(n(T)-5+2)+1>\frac{1}{6}(n(T)+2)$.
If $u$ has a child $v'$ that is distinct from $v$ such that $T_{v'}$ has depth $2$,
then, by symmetry, we may assume that $v'$ has degree $2$ in $T$.
Let $w'$ be the child of $v'$ in $T$.
By Lemma \ref{lemmad31}, some minimum exponential dominating set $S$ of $T$ contains $w$ and $w'$.
If $T'$ arises from $T$ by removing $v'$ together with the at most five descendants of $w$ and $v'$, and adding a new child at $v$, 
then $(S\setminus \{ w,w'\})\cup \{ v\}$ is an exponential dominating set of $T'$, and we obtain, by induction, 
$\gamma_e(T)\geq \gamma_e(T')+1\geq \frac{1}{6}(n(T')+2)+1\geq \frac{1}{6}(n(T)-5+2)+1>\frac{1}{6}(n(T)+2)$.
Hence, we may assume that $v$ is the only child of $u$.

If $u$ has no parent, then $n(T)\leq 5$ and $\gamma_e(T)=2\geq \frac{1}{6}(n(T)+2)$.
Hence, we may assume that $u$ has a parent $x$.
If $u$ is the only child of $x$, then 
either $x$ has no parent, which implies $n(T)\leq 6$ and $\gamma_e(T)=2\geq \frac{1}{6}(n(T)+2)$,
or $x$ has a parent $y$.
In the latter case, Lemma \ref{lemmad31} implies that $T$ has a minimum exponential dominating set $S$ such that $S\cap V(T_x)=\{ w\}$.
Since, by Lemma \ref{lemmad31}, 
$S\setminus \{ w\}$ is an exponential dominating set of $T'=T-V(T_u)$, we obtain, by induction, 
$\gamma_e(T)\geq \gamma_e(T')+1\geq \frac{1}{6}(n(T')+2)+1\geq \frac{1}{6}(n(T)-5+2)+1>\frac{1}{6}(n(T)+2)$.
Hence, we may assume that $x$ has a child $u'$ that is distinct from $u$.
If $u'$ is a leaf, then, by Lemma \ref{lemmad31}, some minimum exponential dominating set $S$ of $T$ contains $x$ and $w$.
Since $S\setminus \{ w\}$ is an exponential dominating set of $T-V(T_u)$, we can argue as above.
Hence, we may assume that $u'$ is not a leaf.
If $T_{u'}$ has depth $1$, 
then, by Lemma \ref{lemmad31}, some minimum exponential dominating set $S$ of $T$ contains $u'$ and $w$.
If $T'$ arises from $T$ by removing the at most six descendants of $u$ and $u'$,
then $(S\setminus \{ u',w\})\cup \{ x\}$ is an exponential dominating set of $T'$, and we obtain, by induction, 
$\gamma_e(T)\geq \gamma_e(T')+1\geq \frac{1}{6}(n(T')+2)+1\geq \frac{1}{6}(n(T)-6+2)+1=\frac{1}{6}(n(T)+2)$.
Hence, we may assume that $T_{u'}$ has depth at least $2$.
If $T_{u'}$ has depth $2$, then, by previous arguments, we may assume that $u'$ has degree $2$ in $T$. 
Furthermore, by Lemma \ref{lemmad31}, some minimum exponential dominating set $S$ of $T$ contains $w$ and the child $v'$ of $u'$.
If $T'$ arises from $T$ by removing the at most seven vertices in $V(T_{u'})\cup V(T_w)$, and adding a new child at $u$,
then $(S\setminus \{ v',w\})\cup \{ u\}$ is an exponential dominating set of $T'$, and we obtain, by induction, 
$\gamma_e(T)\geq \gamma_e(T')+1\geq \frac{1}{6}(n(T')+2)+1\geq \frac{1}{6}(n(T)-6+2)+1=\frac{1}{6}(n(T)+2)$.
Hence, we may assume that $T_{u'}$ has depth at least $3$.
By symmetry, we may assume that $u'$ as a unique child $v'$, which has a unique child $w'$.
By Lemma \ref{lemmad31}, some minimum exponential dominating set $S$ of $T$ contains $w$ and $w'$.
If $T'$ arises from $T$ by removing the at most eight vertices in $V(T_{u'})\cup V(T_w)$, and adding two new children at $v$,
then $(S\setminus \{ w,w'\})\cup \{ v\}$ is an exponential dominating set of $T'$, and we obtain, by induction, 
$\gamma_e(T)\geq \gamma_e(T')+1\geq \frac{1}{6}(n(T')+2)+1\geq \frac{1}{6}(n(T)-6+2)+1=\frac{1}{6}(n(T)+2)$,
which completes the proof. $\Box$

\medskip

\begin{figure}[H]
\begin{center}
\unitlength 0.7mm 
\linethickness{0.4pt}
\ifx\plotpoint\undefined\newsavebox{\plotpoint}\fi 
\begin{picture}(186,21)(0,0)
\put(185,10){\circle*{2}}
\put(5,10){\circle*{2}}
\put(5,20){\circle*{2}}
\put(185,20){\circle*{2}}
\put(15,15){\circle*{2}}
\put(25,15){\circle*{2}}
\put(65,15){\circle*{2}}
\put(105,15){\circle*{2}}
\put(145,15){\circle*{2}}
\put(35,15){\circle*{2}}
\put(75,15){\circle*{2}}
\put(115,15){\circle*{2}}
\put(155,15){\circle*{2}}
\put(45,15){\circle*{2}}
\put(85,15){\circle*{2}}
\put(125,15){\circle*{2}}
\put(165,15){\circle*{2}}
\put(5,20){\line(2,-1){10}}
\put(15,15){\line(-2,-1){10}}
\put(55,15){\circle*{2}}
\put(95,15){\circle*{2}}
\put(135,15){\circle*{2}}
\put(175,15){\circle*{2}}
\put(55,5){\circle*{2}}
\put(95,5){\circle*{2}}
\put(135,5){\circle*{2}}
\put(15,15){\line(1,0){40}}
\put(55,15){\line(1,0){40}}
\put(95,15){\line(1,0){40}}
\put(55,15){\line(0,-1){10}}
\put(95,15){\line(0,-1){10}}
\put(135,15){\line(0,-1){10}}
\put(175,15){\circle*{2}}
\put(185,20){\line(-2,-1){10}}
\put(175,15){\line(2,-1){10}}
\put(145,15){\line(-1,0){10}}
\put(175,15){\line(-1,0){30}}
\end{picture}
\end{center}
\caption{A tree $T$ with $\gamma_e(T)=\frac{1}{5}(n(T)+1)$.}\label{fig2}
\end{figure}
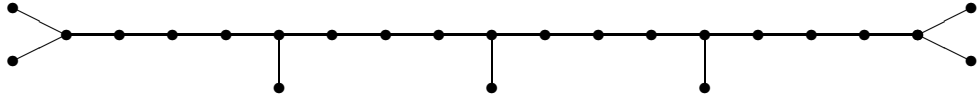

\noindent We believe that the bound in Theorem \ref{theoremlb3} can be improved to 
$\gamma_e(T)\geq \frac{1}{5}(n(T)+1)$.
In view of trees as the one illustrated in Figure \ref{fig2}, 
this bound would be a tight.

\medskip

\noindent There is no lower bound on $\gamma_e(T)$ that is linear in $n(T)$ 
for trees $T$ whose maximum degree is allowed to be $5$ or bigger.

In fact, let $d$ and $h$ be positive integers such that 
\begin{eqnarray}\label{e2}
d&\geq &\frac{4}{3}\cdot 2^{2h-1}+h-1.
\end{eqnarray}
Let $T$ be the rooted tree $T(5,4,\ldots,4)$ of maximum degree $5$ and depth $d$, 
and let $S$ be a set of $5\cdot 4^{d-h-1}$ leaves of $T$ 
that contains exactly one descendant of every vertex of $T$ of depth $d-h$.

If $u$ is a vertex of $T$ of depth at least $d-h$, then
\begin{eqnarray*}
w_{(T,S)}(u) & \geq & 
\left(\frac{1}{2}\right)^{2h-1}\\
&&+
\underbrace{3\cdot 4^0\cdot\left(\frac{1}{2}\right)^{2h+2-1}
+3\cdot 4^1\cdot\left(\frac{1}{2}\right)^{2h+4-1}
+\cdots
+3\cdot 4^{d-h-2}\cdot\left(\frac{1}{2}\right)^{2h+2+2(d-h-2)-1}}_{\mbox{\footnotesize $d-h-1$ terms}}\\
&&+4\cdot 4^{d-h-1}\cdot\left(\frac{1}{2}\right)^{2h+2+2(d-h-1)-1}\\
& = & 
\left(\frac{1}{2}\right)^{2h-1}
+
\underbrace{\frac{3}{4}\cdot \left(\frac{1}{2}\right)^{2h-1}
+\frac{3}{4}\cdot \left(\frac{1}{2}\right)^{2h-1}
+\cdots
+\frac{3}{4}\cdot \left(\frac{1}{2}\right)^{2h-1}}_{\mbox{\footnotesize $d-h-1$ terms}}
+\left(\frac{1}{2}\right)^{2h-1}\\
& > & \frac{3}{4}(d-h+1)\left(\frac{1}{2}\right)^{2h-1}\\
& \geq & 1.
\end{eqnarray*}
Similarly, if $u$ is a vertex of $T$ of depth $i$ for some integer $i$ with $1\leq i\leq d-h$, then
\begin{eqnarray*}
w_{(T,S)}(u) & = & 
4^{d-h-i}\cdot \left(\frac{1}{2}\right)^{d-i-1}\\
&& +
\underbrace{3\cdot 4^0\cdot 4^{d-h-i}\cdot \left(\frac{1}{2}\right)^{d-i+2-1}
+\cdots
+3\cdot 4^{i-2}\cdot 4^{d-h-i}\cdot \left(\frac{1}{2}\right)^{d-i+2+2(i-2)-1}}_{\mbox{\footnotesize $i-1$ terms}}\\
&& +4\cdot 4^{i-1}\cdot 4^{d-h-i}\cdot \left(\frac{1}{2}\right)^{d-i+2+2(i-1)-1}\\
& \geq & \frac{3}{4}\cdot (i+1)\cdot 4^{d-h-i}\cdot \left(\frac{1}{2}\right)^{d-i-1}\\
& \geq & \frac{3}{4}\cdot (d-h+1)\cdot 4^{d-h-d+h}\cdot \left(\frac{1}{2}\right)^{d-d+h-1}\\
& = & \frac{3}{4}\cdot (d-h+1)\cdot \left(\frac{1}{2}\right)^{h-1}\\
& \geq & 1.
\end{eqnarray*}
Finally, $w_{(T,S)}(r)\geq 1$ for the root $r$ of $T$.
Now, selecting, for some large integer $h$, the smallest integer $d$ that satisfies (\ref{e2})
yields a tree $T$ for which $\gamma_e(T)=O\left(\frac{n(T)}{\log(n(T))}\right)$.

We do not know whether or not $\gamma_e(T)=\Omega(n(T))$ for trees $T$ of maximum degree at most $4$.
Furthermore, it may be true that for any fixed positive integer $\Delta$,
we have $\gamma_e(T)=\Omega\left(\frac{n(T)}{\log(n(T))}\right)$ for trees $T$ of maximum degree at most $\Delta$.

\medskip

\noindent {\it Proof of Theorem \ref{theoremlb3gen}:} 
Let $S$ be a minimum exponential dominating set of $G$.
Let $k=|S|$.
For a vertex $u$ in $V(G)\setminus S$, let 
$S(u)=\{ v\in S:{\rm dist}_G(u,v)\leq \log_2(k)+2\}$, and
$$\tilde{w}_{(G,S)}(u)=\sum_{v\in S(u)}\left(\frac{1}{2}\right)^{{\rm dist}_G(u,v)-1}.$$
Note that
\begin{eqnarray*}
\tilde{w}_{(G,S)}(u) 
& = & \left(\sum_{v\in S}\left(\frac{1}{2}\right)^{{\rm dist}_G(u,v)-1}\right)
-\left(\sum_{v\in S\setminus S(u)}\left(\frac{1}{2}\right)^{{\rm dist}_G(u,v)-1}\right)\\
& \geq & \left(\sum_{v\in S}\left(\frac{1}{2}\right)^{{\rm dist}_{(G,S)}(u,v)-1}\right)
-k\left(\frac{1}{2}\right)^{(\log_2(k)+2)-1}\\
& = & w_{(G,S)}(u)-\frac{1}{2}\\
& \geq & \frac{1}{2},
\end{eqnarray*}
and hence,
$$\sum_{u\in V(G)\setminus S}\tilde{w}_{(G,S)}(u)\geq \frac{1}{2}(n(G)-k).$$
Since $G$ has maximum degree at most $3$,
for every vertex $v$ of $G$, 
there are at most $3\cdot 2^{i-1}$ vertices $u$ of $G$ with ${\rm dist}_G(u,v)=i$.
Therefore,
\begin{eqnarray*}
\sum_{u\in V(G)\setminus S}\tilde{w}_{(G,S)}(u) 
& = & \sum_{u\in V(G)\setminus S}\,\,\,\sum_{v\in S(u)}\,\,\,\left(\frac{1}{2}\right)^{{\rm dist}_G(u,v)-1}\\
& = & \sum_{v\in S}\,\,\,\,\,\,\sum_{u\in V(G)\setminus S:\,\, v\in S(u)}\,\,\,\left(\frac{1}{2}\right)^{{\rm dist}_G(u,v)-1}\\
& = & \sum_{v\in S}\,\,\,\sum_{i=1}^{\left\lfloor\log_2(k)+2\right\rfloor}\,\,\,\sum_{u\in V(G)\setminus S:\,\, {\rm dist}_G(u,v)=i}\,\,\,\left(\frac{1}{2}\right)^{i-1}\\
& \leq & \sum_{v\in S}\,\,\,\sum_{i=1}^{\left\lfloor\log_2(k)+2\right\rfloor}\,\,\,\left(\frac{1}{2}\right)^{i-1}\cdot 3\cdot 2^{i-1}\\
& = & \sum_{v\in S}\,\,\,\sum_{i=1}^{\left\lfloor\log_2(k)+2\right\rfloor}\,\,\, 3\\
& = & \sum_{v\in S}3\left\lfloor\log_2(k)+2\right\rfloor\\
& \leq & 3(\log_2(k)+2)k.
\end{eqnarray*}
This implies $\frac{1}{2}(n(G)-k)\leq 3(\log_2(k)+2)k$, and hence, 
$n(G)\leq \left(6\log_2(k)+13\right)k$.
By Theorem \ref{theoremub3}, we have $6\log_2(k)+13\leq 6\log_2(n(G)+2)+4$, and hence,
$k\geq \frac{n(G)}{6\log_2(n(G)+2)+4}$, which completes the proof. $\Box$

\medskip

\noindent {\it Proof of Theorem \ref{theoremgirth}:}
For a positive integer $d$, 
let $T$ be the rooted tree in which every leaf has depth $d$, and every vertex that is not a leaf, has degree $3$.
Let $S$ be a random subset of $V(T)$ 
that contains each vertex of $T$ independently at random with probability $p$ for some $p\in [0,1]$.

Let $r$ be the root of $T$.
If $r\in S$, then $w_{(T,S)}(r)=2$.

Let $v$ be a vertex at distance $i$ from $r$ for some $i\in [d]$.
If $v$ is the only vertex on the path of order $i+1$ in $T$ between $r$ and $v$ that belongs to $S$,
then $v$ contributes $\left(\frac{1}{2}\right)^{i-1}$ to $w_{(T,S)}(r)$;
otherwise $v$ contributes $0$ to $w_{(T,S)}(r)$.
Therefore, the vertex $v$ contributes $\left(\frac{1}{2}\right)^{i-1}$ 
to $w_{(T,S)}(r)$ exactly with probability $p\cdot (1-p)^i$.
Since there are exactly $3\cdot 2^{i-1}$ vertices at distance $i$ from $r$ for every $i\in [d]$,
we obtain, by linearity of expectation,
\begin{eqnarray*}
{\bf E}\left[w_{(T,S)}(r)\right]
&=& 2p+\sum_{i=1}^d 3\cdot 2^{i-1}\cdot \left(\frac{1}{2}\right)^{i-1}\cdot p\cdot (1-p)^i\\
&=& 2p+3(1-p)\left(1-(1-p)^d\right)\\
& \geq & 2p+3(1-p)\left(1-e^{-pd}\right)\\
& \geq & 3-p-3e^{-pd}.
\end{eqnarray*}
By Lemma \ref{lemmad31}, we have $w_{(T,S)}(r)\leq 3$, and hence,
\begin{eqnarray*}
{\bf E}\left[w_{(T,S)}(r)\right]
&\leq &{\bf P}\left[w_{(T,S)}(r)<1\right]+3{\bf P}\left[w_{(T,S)}(r)\geq 1\right]\\
&=& 3-2{\bf P}\left[w_{(T,S)}(r)<1\right],
\end{eqnarray*}
which implies 
$${\bf P}\left[w_{(T,S)}(r)<1\right]\leq \frac{p}{2}+\frac{3}{2}e^{-pd}.$$
Now, let $G$ be a cubic graph of girth at least $2d+1$.
Let $S_0$ be a random subset of $V(G)$ that contains each vertex of $G$ independently at random with probability $p$.
Let $S_1=\{ u\in V(G):w_{(G,S_0)}(u)<1\}.$
By Lemma \ref{lemmad31}, $S_0\cup S_1$ is an exponential dominating set of $G$.

Let $u$ be a vertex of $G$.
By the girth condition, the subgraph of $G$ induced by the vertices at distance at most $d$ from $u$ is isomorphic to $T$. Let $S$ be the restriction of $S_0$
to these vertices. We have
$${\bf P}\left[w_{(G,S_0)}(u)<1\right]\leq {\bf P}\left[w_{(T,S)}(r)<1\right]\leq \frac{p}{2}+\frac{3}{2}e^{-pd},$$
and, by linearity of expectation and the first moment method, we obtain 
\begin{eqnarray*}
\gamma_e(G)
&\leq & {\bf E}\left[|S_0|\right]+{\bf E}\left[|S_1|\right]\\
&\leq & p n(G)+\left(\frac{p}{2}+\frac{3}{2}e^{-pd}\right)n(G)\\
&= & \frac{3}{2}\left(p+e^{-pd}\right)n(G).
\end{eqnarray*}
Now, let $\epsilon$ be such that $0<\epsilon<1$.
For $p(\epsilon)=\frac{\epsilon}{3}$ and $d(\epsilon)=\left\lceil\frac{3}{\epsilon}\ln\left(\frac{3}{\epsilon}\right)\right\rceil$, 
we obtain that $\frac{3}{2}\left(p(\epsilon)+e^{-p(\epsilon)d(\epsilon)}\right)\leq \epsilon$.
Therefore, every cubic graph $G$ of girth at least $2d(\epsilon)+1$
satisfies $\gamma_e(G)\leq \epsilon n(G)$.

Finally, let $\alpha$ be such that $0<\alpha<\frac{2}{3\ln(2)}$.
It is known \cite{bh,w} that there are arbitrarily large cubic graphs $G$ of girth $g(G)\geq \frac{4}{3\ln(2)}\ln(n(G))-2$.
Since $\alpha<\frac{2}{3\ln(2)}\approx 0.96$, there are infinitely many such graphs $G$ 
with
$\ln(\ln(n(G)))\leq \ln(n(G))^{1-\alpha}$
and
$\lceil \alpha \ln(n(G))\rceil\leq \frac{2}{3\ln(2)}\ln(n(G))-\frac{3}{2}$.
For $p(\alpha)=\frac{\ln(\ln(n(G)))}{\ln(n(G))}$ and 
$d(\alpha)=\lceil \alpha \ln(n(G))\rceil$, 
we obtain
$p(\alpha)\leq \frac{1}{\ln(n(G))^{\alpha}}$
and
$\frac{g(G)-1}{2}\geq \frac{2}{3\ln(2)}\ln(n(G))-\frac{3}{2}\geq d(\alpha)$.
This implies
$\frac{3}{2}\left(p(\alpha)+e^{-p(\alpha)d(\alpha)}\right)\leq \frac{3}{\ln(n(G))^{\alpha}}$,
and hence
$\gamma_e(G)\leq \frac{3n(G)}{\ln(n(G))^{\alpha}}$.
$\Box$

\medskip

\noindent {\it Proof of Theorem \ref{theoremlb3sp}:}
Clearly, we may assume that $S$ does not contain all vertices of $G$.
Let $H$ arise from $G$ by removing all edges between vertices in $S$, 
and replacing each vertex $u$ in $S$ that has $d$ neighbors in $V(G)\setminus S$ by $d$ copies of degree $1$,
each copy being adjacent to one of the neighbors of $u$ in $V(G)\setminus S$.
Let $\bar{S}=V(H)\setminus (V(G)\setminus S)$.
Clearly, $w_{(H,\bar{S})}(u)=w_{(G,S)}(u)\geq 3$ for every vertex $u\in V(H)\setminus \bar{S}$.

Considering the at most three neighbors of each vertex $v$ in $V(H)\setminus \bar{S}$, 
Lemma \ref{lemmad31} implies that all components of $H$ are trees 
whose internal vertices all have degree $3$ and belong to $V(H)\setminus \bar{S}$,
and whose endvertices belong to $\bar{S}$.
Furthermore, no two endvertices in one component of $H$ are copies of the same vertex in $S$.
Note that every tree with vertices of degree $3$ and $1$ only, 
has order $2\ell-2$, where $\ell$ is the number of its endvertices.
Let $H$ have $k$ components.
Since $H$ has $|\bar{S}|$ endvertices and $n(G)-|S|$ internal vertices, 
we obtain $|\bar{S}|=n(G)-|S|+2k$.

If some vertex in $S$ has $3$ neighbors in $V(G)\setminus S$, then $k\geq 3$.
Since $|\bar{S}|\leq 3|S|$, this implies 
$3|S|\geq n(G)-|S|+6$, and hence, $|S|\geq \frac{1}{4}(n(G)+6)$.
If no vertex in $S$ has $3$ neighbors in $V(G)\setminus S$, 
then, 
either $k\geq 2$ and $|\bar{S}|\leq 2|S|$,
or $k\geq 1$ and $|\bar{S}|\leq |S|$.
In both cases, this implies $|S|>\frac{1}{4}(n(G)+6)$.
$\Box$

\medskip

\noindent It is not difficult to see that the extremal graphs $G$ for Theorem \ref{theoremlb3sp}
arise from three disjoint trees $T_1$, $T_2$, and $T_3$ with vertices of degree $3$ and $1$ only, 
and $\frac{1}{4}(n(G)+6)$ endvertices each,
by identifying each leaf of $T_1$ with one leaf of $T_2$ and one leaf of $T_3$.
For instance $K_{3,3}$ is one such graph.

\medskip

\noindent We proceed to the proof of Theorem \ref{theoremalg3}.
For a given tree $T$ of maximum degree at most $3$, 
and a given set $S$ of vertices of $T$ such that 
\begin{itemize}
\item $S$ is not an exponential dominating set of $T$, but 
\item $S$ is a subset of some minimum exponential dominating set of $T$,
\end{itemize}
we explain how to identify, in polynomial time, a vertex $u$ in $V(T)\setminus S$ 
such that $S\cup \{ u\}$ is a subset of some minimum exponential dominating set of $T$.
Iteratively applying this extension starting with the empty set, 
clearly yields a proof of Theorem \ref{theoremalg3}.

Since the different maximal subtrees $T'$ of $T$ 
for which $S$ contains only endvertices of $T'$ can be handled completely independently,
we may assume that $S$ is a subset of the set of endvertices of $T$.
We root $T$ in some vertex $r$ in $V(T)\setminus S$.

For every vertex $u$ of $T$, let
$$\partial w(u)=\max\left\{ 2^{{\rm dist}_T(u,v)}\left(1-w_{(T_u,S\cap V(T_u))}(v)\right): v\in V(T_u)\right\}.$$
Clearly, given $T$ and $S$, the value of $\partial w(u)$ can be calculated efficiently for every vertex $u$.

\begin{lemma}\label{lemmaalg3}
Let $T$, $S$, and $\partial w$ be as above.
\begin{enumerate}[(i)]
\item If $u$ is a vertex of $T$ that is distinct from the root $r$
such that $\partial w(u)>1$ and $\partial w(v)\leq 1$ for every descendant $v$ of $u$ in $T$, 
then $S\cup \{ u\}$ is a subset of some minimum exponential dominating set of $T$.
\item If $\partial w(u)\leq 1$ for every vertex $u$ of $T$, 
then $S\cup \{ r\}$ is a minimum exponential dominating set of $T$.
\end{enumerate}
\end{lemma}
{\it Proof:} Let $\bar{S}$ be a minimum exponential dominating set of $T$ such that $S\subseteq \bar{S}$.

\medskip

\noindent (i) Our first goal is to show that $\bar{S}\setminus S$ intersects $V(T_u)$.
Therefore, suppose that $\bar{S}\cap V(T_u)=S\cap V(T_u)$.
Let $v\in V(T_u)$ be such that $\partial w(u)=2^{{\rm dist}_T(u,v)}\left(1-w_{(T_u,S\cap V(T_u))}(v)\right)>1$,
which implies $w_{(T_u,S\cap V(T_u))}(v)<1-\left(\frac{1}{2}\right)^{{\rm dist}_T(u,v)}$.
Let $u^-$ be the parent of $u$ in $T$.
By Lemma \ref{lemmad31}, $w_{(T-uu^-,\bar{S})}(u^-)\leq 2$.
Since $\bar{S}\cap V(T_u)=S\cap V(T_u)$ and, consequently, $\bar{S}$ does not intersect the path in $T$ between $u^-$ and $v$,
we obtain
\begin{eqnarray*}
w_{(T,\bar{S})}(v) & = & w_{(T-uu^-,\bar{S})}(v)+\left(\frac{1}{2}\right)^{{\rm dist}_T(u^-,v)}w_{(T-uu^-,\bar{S})}(u^-)\\
& = & w_{(T_u,S\cap V(T_u))}(v)+\left(\frac{1}{2}\right)^{{\rm dist}_T(u,v)+1}w_{(T-uu^-,\bar{S})}(u^-)\\
& < & \left(1-\left(\frac{1}{2}\right)^{{\rm dist}_T(u,v)}\right)+2\cdot \left(\frac{1}{2}\right)^{{\rm dist}_T(u,v)+1}\\
&=& 1,
\end{eqnarray*}
which is a contradiction.
Hence, $\bar{S}$ intersects $V(T_u)\setminus S$.

Let $S'=S\cup \{ u\}\cup (\bar{S}\setminus V(T_u))$.
Clearly, $|S'|\leq |\bar{S}|$.

Our next goal is to show that $S'$ is an exponential dominating set of $T$.
By Lemma \ref{lemmad31}, $w_{(T-uu^-,\bar{S})}(u)\leq 2$.
Since $w_{(T-uu^-,S')}(u)=2$, 
this implies $w_{(T,S')}(v)\geq w_{(T,\bar{S})}(v)\geq 1$ for every vertex $v$ in $V(T)\setminus V(T_u)$.
Trivially, $w_{(T,S')}(u)\geq 1$.
Now, let $v\in V(T_u)\setminus \{ u\}$.
Let $u^+$ be the child of $u$ on the path in $T$ between $u$ and $v$.
Since $v\in V(T_{u^+})$, we have
$$1\geq \partial w(u^+)\geq 2^{{\rm dist}_T(u^+,v)}\left(1-w_{(T_{u^+},S\cap V(T_{u^+}))}(v)\right),$$
which implies
\begin{eqnarray*}
w_{(T,S')}(v) & = & w_{(T_{u^+},S\cap V(T_{u^+}))}(v)+\left(\frac{1}{2}\right)^{{\rm dist}_T(u,v)-1}\\
& \geq & 1-\left(\frac{1}{2}\right)^{{\rm dist}_T(u^+,v)}+\left(\frac{1}{2}\right)^{{\rm dist}_T(u,v)-1}\\
& = & 1.
\end{eqnarray*}
Hence $S'$ is an exponential dominating set of $T$.

Since $|S'|\leq |\bar{S}|$, this completes the proof of (i).

\medskip

\noindent (ii) As in the proof of (i), it follows that $S\cup \{ r\}$ is an exponential dominating set of $T$.
By hypothesis, $S$ is not an exponential dominating set of $T$, which implies $|S|<|\bar{S}|$.
Hence $S\cup \{ r\}$ is a minimum exponential dominating set of $T$.
$\Box$

\medskip

\noindent {\it Proof of Theorem \ref{theoremnpc}:}
In order to prove APX-hardness, we use that {\sc Min Vertex Cover} is APX-complete for cubic graphs \cite{ak}.
For a given cubic graph $G$, we construct a subcubic graph $H$ such that 
$n(H)=16n(G)$ and $\gamma_e(H)=\tau(G)+3n(G)$, 
where $\tau(G)$ is the minimum order of a vertex cover of $G$.
Furthermore, given an exponential dominating set $S$ of $H$,
we explain how to construct efficiently a vertex cover $C$ of $G$ with $|C|\leq |S|-3n(G)$.
Note that since $G$ is cubic, we have $m(G)=\frac{3}{2}n(G)$, which implies $\tau(G)\geq \frac{n(G)}{2}$,
and hence, $\gamma_e(H)\geq \frac{7}{2}n(G)=\frac{7}{32}n(H)$.

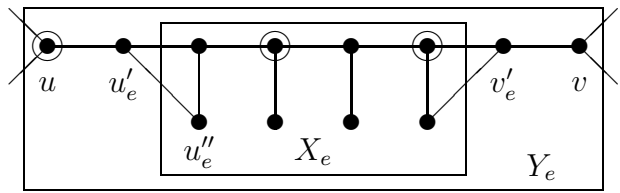
\begin{figure}[H]
\begin{center}
\unitlength 1mm 
\linethickness{0.4pt}
\ifx\plotpoint\undefined\newsavebox{\plotpoint}\fi 
\begin{picture}(80,24)(0,0)
\put(5,17){\circle*{2}}
\put(5,17){\circle{4}}
\put(15,17){\circle*{2}}
\put(25,17){\circle*{2}}
\put(35,17){\circle*{2}}
\put(35,17){\circle{4}}
\put(45,17){\circle*{2}}
\put(55,17){\circle*{2}}
\put(55,17){\circle{4}}
\put(65,17){\circle*{2}}
\put(75,17){\circle*{2}}
\put(35,7){\circle*{2}}
\put(45,7){\circle*{2}}
\put(5,17){\line(1,0){70}}
\put(75,17){\line(1,1){5}}
\put(75,17){\line(1,-1){5}}
\put(45,17){\line(0,-1){10}}
\put(35,17){\line(0,-1){11}}
\put(5,17){\line(-1,1){5}}
\put(5,17){\line(-1,-1){5}}
\put(25,7){\circle*{2}}
\put(55,7){\circle*{2}}
\put(25,17){\line(0,-1){10}}
\put(25,7){\line(-1,1){10}}
\put(55,17){\line(0,-1){10}}
\put(55,7){\line(1,1){10}}
\put(20,0){\framebox(40,20)[cc]{}}
\put(2,-2){\framebox(76,24)[cc]{}}
\put(5,12){\makebox(0,0)[cc]{$u$}}
\put(15,12){\makebox(0,0)[cc]{$u'_e$}}
\put(25,3){\makebox(0,0)[cc]{$u''_e$}}
\put(75,12){\makebox(0,0)[cc]{$v$}}
\put(65,12){\makebox(0,0)[cc]{$v'_e$}}
\put(40,3){\makebox(0,0)[cc]{$X_e$}}
\put(70,1){\makebox(0,0)[cc]{$Y_e$}}
\end{picture}
\end{center}
\caption{The gadget $G_e$ for the edge $e=uv$.}\label{fig1}
\end{figure}

Let $G$ be a cubic graph.
Let $H$ arise from $G$ by replacing every edge $e=uv$ of $G$ by a copy $G_e$ of the subgraph shown in Figure \ref{fig1}.
Clearly, the order of $H$ is $n(G)+10m(G)=16n(G)$.

\noindent Let $X$ be a vertex cover of $G$.
For every edge $e$, arbitrarily select one vertex $x(e)$ in $X$ that is incident with $e$.
Starting with $S=\emptyset$, add to $S$ all vertices of $X$.
Furthermore, for every edge $e$ of $G$, if $e=uv$ and $x(e)=u$, 
then add to $S$ the two vertices from the set $X_e$ indicated in Figure \ref{fig1}.
Since $G$ is cubic, and $X$ is a vertex cover of $G$, the resulting set $S$ is an exponential dominating set of $G$,
and hence, $\gamma_e(H)\leq\tau(G)+2m(G)=\tau(G)+3n(G)$.

\medskip

\noindent Let $S$ be an exponential dominating set of $H$.
We may assume that there is no exponential dominating set $S'$ of $H$ with $|S'|\leq |S|$ and $|S\setminus S'|+|S'\setminus S|\leq 24$
such that $|S'\cap V(G)|>|S\cap V(G)|$.
Otherwise, such a set $S'$ can be found efficiently, and iteratively replacing $S$ with $S'$
yields an exponential dominating set with the desired property after at most $n(G)$ steps.

Let $e=uv$ be an edge of $G$.
By Lemma \ref{lemmad31}, 
the set $S$ contains at least two of the eight vertices of the set $X_e$ indicated in Figure \ref{fig1}.

If $S$ does not contain $u$ but contains the neighbor $u'_e$ of $u$ in $G_e$, 
then replacing $u'_e$ and the (at least two) vertices in $S\cap X_e$
with the three vertices indicated in Figure \ref{fig1} 
yields a minimum exponential dominating set $S'$ of $H$ with $|S'\cap V(G)|>|S\cap V(G)|$,
which is a contradiction.
If $S$ does not contain $u$ but contains the two neighbors of $u'_e$ in $X_e$,
then replacing the neighbor $u''_e$ of $u'_e$ in $G_e$ with $u$ 
yields a similar contradiction.
Hence, if $S$ does not contain $u$, then $S$ does not contain $u'_e$,
and $S$ does not contain both neighbors of $u'_e$ in $X_e$.

If $S$ contains only two of the twelve vertices of the set $Y_e$ indicated in Figure \ref{fig1},
then $S\cap (Y_e\setminus X_e)$ is empty, 
and we may assume, by symmetry, that $S$ contains neither $u'_e$ nor any neighbor of $u'_e$.
Since $w_{(H,S)}(u''_e)\geq 1$, this implies $w_{(H-uu'_e,S)}(u)\geq 2$.
Applying Lemma \ref{lemmad31} to the vertex $u$ of degree $2$ in the subcubic graph $H-uu'_e$ 
implies that there is a neighbor $w$ of $u$ distinct from $v$ 
such that $S$ either contains $u'_{uw}$
or contains both neighbors of $u'_{uw}$ within $X_{uw}$,
which is a contradiction.
Hence, $S$ contains at least three of the twelve vertices of the set $Y_e$.

If $S$ contains exactly three vertices from $Y_e$ but neither $u$ nor $v$,
then we may assume, by symmetry, that $v'_e$ is not contained in $S$.
Replacing the three vertices in $S\cap Y_e$ with the three vertices indicated in Figure \ref{fig1}
yields a minimum exponential dominating set $S'$ of $H$ with $|S'\cap V(G)|>|S\cap V(G)|$,
which is a contradiction.
Similarly, if $S$ contains at least four vertices from $Y_e$ but neither $u$ nor $v$,
then replacing the vertices in $S\cap Y_e$ with the three vertices indicated in Figure \ref{fig1} and additionally adding $v$
yields a contradiction.
Hence, $S$ contains either $u$ or $v$.
This implies that $S\cap V(G)$ is a vertex cover of $G$ of order at most $|S|-2m(G)=|S|-3n(G)$.
Hence, $\gamma_e(H)\geq\tau(G)+3n(G)$, and the proof is complete. 
$\Box$

\section{Open problems}

We collect the several open problems that are scattered throughout the paper.

Does Conjecture \ref{conjectured3} hold?
More generally, what are the extremal graphs for Theorem \ref{theoremub3}?
Can Theorem \ref{theoremlb3} be improved as explained before Figure \ref{fig2}?
What is the right order of magnitude for the lower bound in Theorem \ref{theoremlb3gen}?
Is $\gamma_e(T)=\Omega(n(T))$ for trees $T$ of maximum degree at most $4$?
Is $\gamma_e(T)=\Omega\left(\frac{n(T)}{\log(n(T))}\right)$ for trees $T$ of bounded maximum degree?
Can Theorem \ref{theoremalg3}, that is, the efficient algorithm, be extended to all trees
or at least to the porous version of exponential domination on subcubic trees?

\end{document}